\documentclass[12pt,a4paper]{article}
\usepackage{latexsym,amsmath,amsthm,amssymb,mathrsfs}
\usepackage[all]{xy}
\usepackage{textcomp}
\usepackage{graphicx}

\begin{document}

\title{Landscape disruption effects in a meta-epidemic model with steady state demographics and migrations saturation.}

\author{Veronica Aimar, Sara Borlengo, Silvia Motto and Ezio Venturino\\
Dipartimento di Matematica ``Giuseppe Peano'',\\
Universit\`a di Torino, \\via Carlo Alberto 10, 10123 Torino, Italy
}
\date{}

\maketitle

\begin{abstract}
We continue the investigations of an ecosystem where a epidemic-affected population can move between two connected patches,
\cite{ABMV}, 
by considering what happens to the system when the migration paths are interrupted
in one direction, or when the infected are not able to exert the effort for migrating into the other patch.
\end{abstract}

\section{Introduction}\label{sec:intro}
In \cite{ABMV} an epidemic-affected one species metapopulation model with fixed size and
immigrations depending inversely on the crowding of the arrival environment has
been introduced, along lines that allow disease consideration in fragmented habitats, \cite{EV}.
Here, we add reproduction capabilities and
specialize the system to two particular cases, when the migrations can occur only in
one direction, or when infected are too weak to undertake any migrating effort.

\section{Unidirectional migrations}
\label{subsec:no_12}
\begin{figure}[!h]
\centering
\includegraphics[scale=0.3]{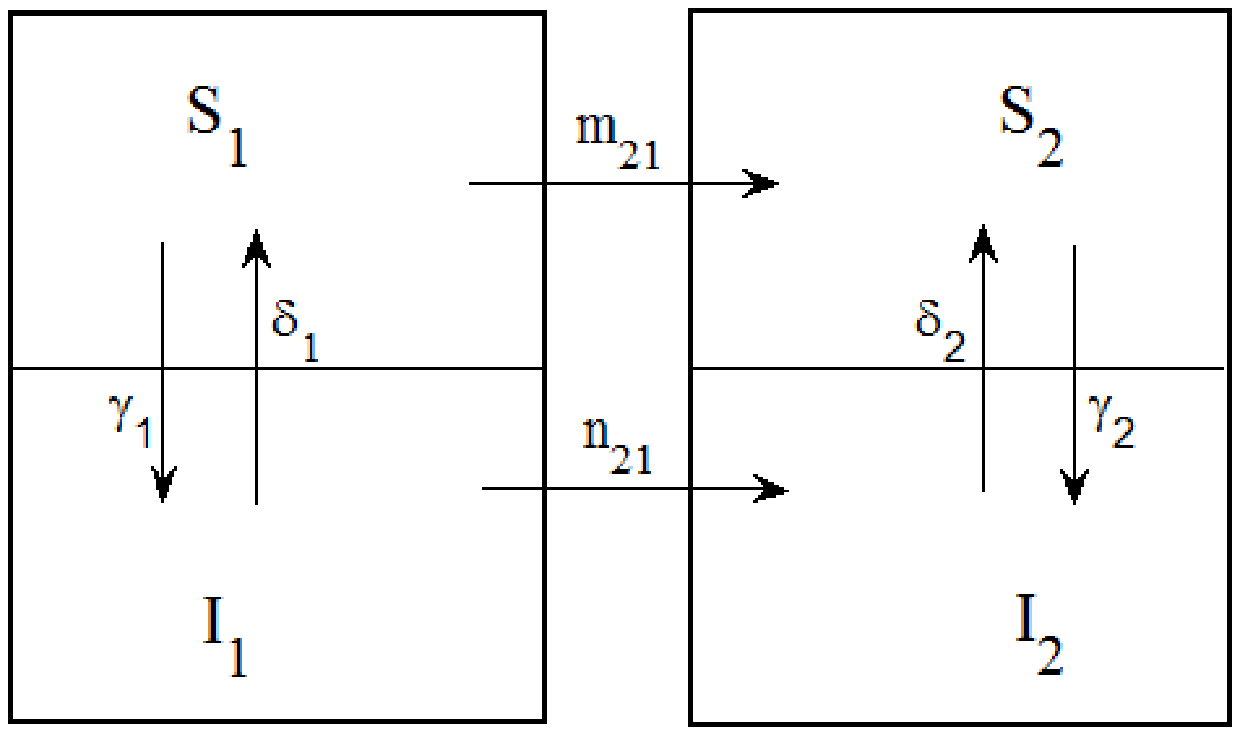} \qquad
\includegraphics[scale=0.3]{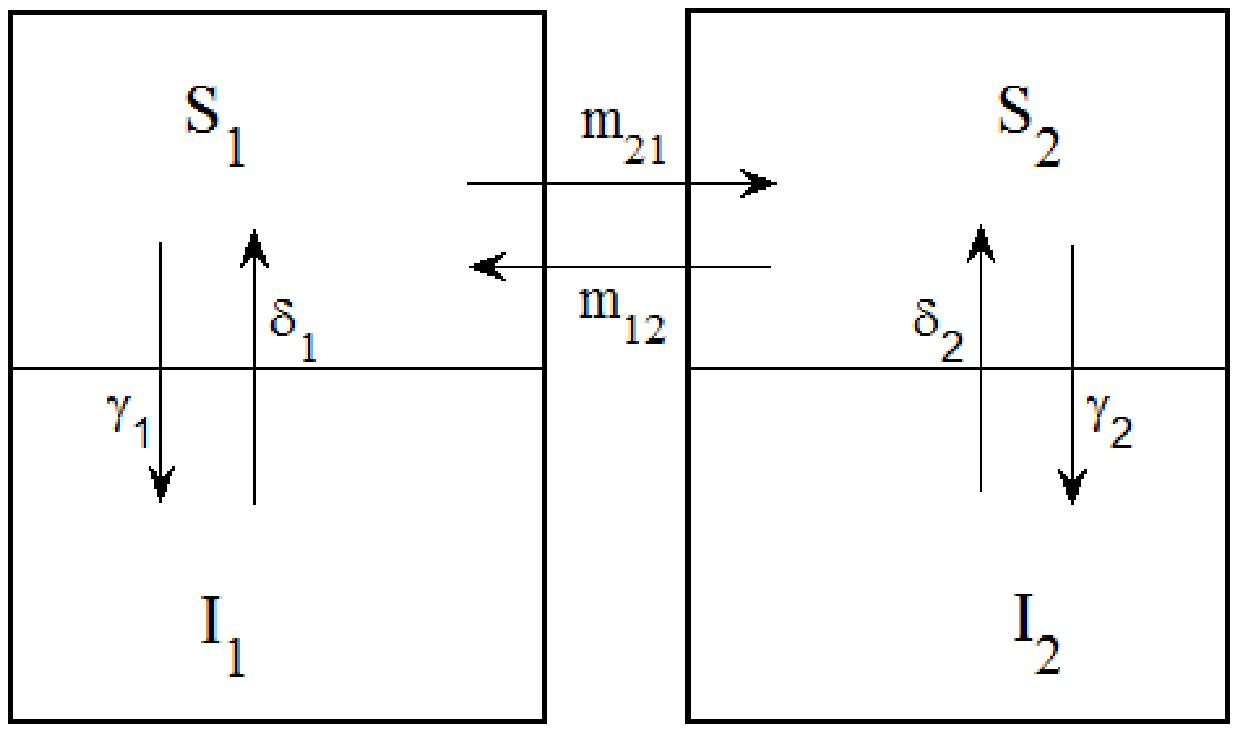}
\caption{Left: no migrations from patch 2 into patch 1.
Right: Infected do not migrate.}
\label{fig:graf_interr_migr_2_1}
\end{figure}

Assume it is not possible to return to patch 1 from the second one.
%; then $m_{12}=n_{12}=0$.
The system is pictured in Figure
\ref{fig:graf_interr_migr_2_1} left.
The model reads
\begin{eqnarray}\label{model_1}
\dot{S_1}=r_1S_1-\gamma_{1}S_1I_1+\delta_1I_1-m_{21}\frac {S_1}{A+I_2+S_2}, \quad
\dot{I_1}=\gamma_{1}S_1I_1-(\delta_1+\mu_1)I_1-n_{21}\frac {I_1}{B+I_2+S_2},\\
\dot{S_2}=r_2S_2-\gamma_{2}S_2I_2+\delta_2I_2+m_{21}\frac {S_1}{A+I_2+S_2}, \quad
\dot{I_2}=\gamma_{2}S_2I_2-(\delta_2+\mu_2)I_2+n_{21}\frac {I_1}{B+I_2+S_2},
\end{eqnarray}
where $r_k$, $k=1,2$ represent the net reproduction rates of the population in each
environment, which is assumed to have different ecological characteristics.
The other parameters have the following meanings
$\mu_k$ is the infected mortality rate in each patch,
$\gamma_k$ the disease contact rate,
$\delta_k$ is the disease recovery rate,
$A$ is the half saturation constant for the susceptibles, and $B$ the one for the infected;
finally the migration rates from patch $j$ into patch $i$ are
$m_{ij}$ for the susceptibles and
$n_{ij}$ for the infected.

The equilibria are the origin, trivially, possibly the coexistence in both patches
with an endemic disease, and the point with only the arrival patch populated by
both susceptibles and infected, $X_1=(0,0,\tilde{S_2},\tilde{I_2})$, 
$$
\tilde{S_2}=\frac{\delta_2+\mu_2}{\gamma_{2}}, \quad
\tilde{I_2}=\frac{r_2(\delta_2+\mu_2)}{\gamma_{2}\mu_2},
$$
which is clearly unconditionally feasible,
and the point $X_2=(\tilde{S_1},0,\tilde{S_2},\tilde{I_2})$, with the first patch
disease-free,
$$
\tilde{S_2}=\frac{\delta_2+\mu_2}{\gamma_2}, \quad
\tilde{I_2}=\frac{\gamma_2m_{21}-r_1\gamma_2A-(\delta_2+\mu_2)r_1}{r_1\gamma_2},
$$
$$
\tilde{S_1}=\frac{-(\delta_2+\mu_2)r_1r_2+\mu_2(\gamma_2m_{21}-r_1\gamma_2A-(\delta_2+\mu_2)r_1)}{r_1\gamma_2}.
$$
This equilibrium is feasible for
\begin{equation}\label{X2_feas}
\gamma_{2}m_{21} \ge r_1\gamma_{2}A+(\delta_2+\mu_2)r_1, \quad
\mu_2[\gamma_{2}m_{21}-r_1\gamma_{2}A-(\delta_2+\mu_2)r_1] \ge (\delta_2+\mu_2)r_1r_2.
\end{equation}

For the equilibrium with both patches populated and with endemic disease, let us sum the
first and third equations of (\ref{model_1}) as well as the second and fourth one, to obtain
\begin{equation}\label{eq2piu4}
r_1\tilde{S_1}-\gamma_{1}\tilde{S_1}\tilde{I_1}+\delta_1\tilde{I_1}+r_2\tilde{S_2}-\gamma_{2}\tilde{S_2}\tilde{I_2}+\delta_2\tilde{I_2}=0, \quad
\gamma_{1}\tilde{S_1}\tilde{I_1}-(\delta_1+\mu_1)\tilde{I_1}+\gamma_{2}\tilde{S_2}\tilde{I_2}-(\delta_2+\mu_2)\tilde{I_2}=0.
\end{equation}
Adding these equations further and solving for $\tilde{S_1}$ as function of $\tilde{I_1}, \tilde{S_2}, \tilde{I_2}$
we substitute it into the second one of (\ref{eq2piu4}) to get
$$
\tilde{S_1}=\frac{-r_2\tilde{S_2}+\mu_1\tilde{I_1}+\mu_2\tilde{I_2}}{r_1}, \quad
\tilde{S_2}=\frac{r_1((\delta_1+\mu_1)\tilde{I_1}+(\delta_2+\mu_2)\tilde{I_2})-\gamma_1\mu_1{\tilde{I_1}}^2-\gamma_1\mu_1\tilde{I_1}\tilde{I_2}}{r_1\gamma_2\tilde{I_2}-r_2\gamma_1\tilde{I_1}}.
$$
Necessary conditions for the feasibility of this equilibrium are either one
of the following two sets of inequalities
\begin{equation}\label{coex_1}
\tilde{I_1}>\frac{r_1\gamma_2\tilde{I_2}}{r_2\gamma_1} , \quad
\tilde{S_2}<\frac{\mu_1\tilde{I_1}+\mu_2\tilde{I_2}}{r_2} , \quad
\tilde{I_2}>\frac{\gamma_1\mu_1{\tilde{I_1}}^2-r_1(\delta_1+\mu_1)\tilde{I_1}}{r_1(\delta_2+\mu_2)-\gamma_1\mu_2\tilde{I_1}} \equiv Z;
\end{equation}
\begin{equation}\label{coex_1bis}
\tilde{I_1}<\frac{r_1\gamma_2\tilde{I_2}}{r_2\gamma_1} , \quad
\tilde{S_2}>\frac{\mu_1\tilde{I_1}+\mu_2\tilde{I_2}}{r_2} , \quad
\tilde{I_2}<Z.
\end{equation}
But we need also to ensure that $Z>0$, so that finally we also get either one of the inequalities
\begin{equation}\label{coex_2}
\frac {\delta_1+\mu_1}{\gamma_1\mu_1} < I_1 < \frac {\delta_2+\mu_2}{\gamma_1\mu_2} ; \quad
\frac {\delta_2+\mu_2}{\gamma_1\mu_2} < I_1 < \frac {\delta_1+\mu_1}{\gamma_1\mu_1} .
\end{equation}

The Jacobian of (\ref{model_1}) is
$$
J=
\left[
\begin{array}{cccc}
-\gamma_1I_1-\eta_1+r_1&-\gamma_1S_1+\delta_1 & \eta_2S_1& \eta_2S_1\\
\gamma_1I_1& \gamma_1S_1-\delta_1-\mu_1-\theta_1& \theta_2I_1& \theta_2I_1\\
\eta_1& 0& -\gamma_2I_2-\eta_2S_1+r_2& -\gamma_2S_2+\delta_2-\eta_2S_1\\
0&\theta_1&\gamma_2I_2-\theta_2I_1& \gamma_2S_2-\delta_2-\mu_2-\theta_2I_1\\
\end{array}
\right]
$$
where
$$
\eta_1=\frac{m_{21}}{A+S_2+I_2}, \quad \eta_2=\frac{m_{21}}{{(A+S_2+I_2)}^2},
\quad
\theta_1=\frac{n_{21}}{B+S_2+I_2}, \quad \theta_2=\frac{n_{21}}{{(B+S_2+I_2)}^2}.
$$

The origin is unstable, since the eigenvalues are 
$r_2$, $-\delta_2-\mu_2$, $(r_1A-m_{21})A^{-1}$,
$-(\delta_1B+n_{21}+\mu_1B ) B^{-1}$.

At $X_1$ we have instead one rather complicated but negative eigenvalue, $\lambda_1<0$ and
$$
\lambda_2=\frac{-\mu_2m_{21}\gamma_2+\mu_2r_1A\gamma_2+\mu_2r_1\delta_2+r_1{\mu_2}^2+r_1r_2(\delta_2+\mu_2)}{\mu_2(A\gamma_2+\delta_2+\mu_2\delta_2+r_2)},
$$
$$
\lambda_{3,4}=\frac{-r_2\delta_2 \pm \sqrt{{r_2}^2{\delta_2}^2-4{\mu_2}^2r_2(\mu_2+\delta_2)}}{2\mu_2}<0.
$$
Stability is then obtained for
\begin{equation}\label{stab_X1}
\mu_2r_1A\gamma_2+\mu_2r_1\delta_2+r_1{\mu_2}^2+r_1r_2(\delta_2+\mu_2)<\mu_2m_{21}\gamma_2.
\end{equation}

At $X_2$, one eigenvalue is explicit,
$\lambda_1=\gamma_1\tilde{S_1}-\delta_1-\mu_1- n_{21}[B+\tilde{S_2}+\tilde{I_2}]^{-1}$
while the remaining ones are the roots of the cubic equation
$\lambda^3+a_2\lambda^2+a_1\lambda+a_0=0$, with
$a_2=\gamma_2\tilde{I_2}+m_{21}\tilde{S_1}D^{-2}-r_2-\gamma_2\tilde{S_2}-\delta_2+\mu_2 $ and
\begin{eqnarray*}
a_1=\left(\frac{m_{21}}{D}-r_1\right), \quad
D=A+\tilde{S_2}+\tilde{I_2} \\
a_2=\gamma_2\mu_2\tilde{I_2}+r_2(\gamma_2\tilde{S_2}-\delta_2-\mu_2)+\frac{m_{21}\tilde{S_1}}{D^2}\left(-\gamma_2\tilde{S_2}
+\delta_2+\mu_2+\gamma_2\tilde{I_2}+\frac{m_{21}}{D}\right),\\
a_0=\left(\frac{m_{21}}{D}-r_1\right)
\left[\gamma_2\tilde{I_2}\left(\mu_2+\frac{m_{21}\tilde{S_1}}{D^2}\right)+\left(r_2-\frac{m_{21}\tilde{S_1}}{D^2}\right)\left(\gamma_2\tilde{S_2}-\delta_2
-\mu_2\right)\right]\\
+\frac{{m_{21}}^2}{D^3}\tilde{S_1}\left(\gamma_2\tilde{S_2}-\gamma_2\tilde{I_2}-\delta_2-\mu_2\right).
\end{eqnarray*}
The Routh-Hurwitz conditions and negativity of the first eigenvalue guarantee stability for
\begin{equation}\label{stab21_caso14}
\lambda_1<0 , \quad
a_0>0 , \quad
a_2>0 , \quad
a_2a_1>a_0 .
\end{equation}
A Hopf bifurcation would be possible if $a_2a_1=a_0$.

The points $X_1$, $X_2$ and coexistence can stably be achieved respectively by the following
parameter choices
\begin{eqnarray*}
r_1=2, \quad r_2=1, \quad \gamma_1=0.5, \quad \gamma_2=1, \quad \delta_1=0.5, \quad \delta_2=2, \quad
\mu_1=\mu_2=1, \\ m_{21}=20, \quad n_{21}=0.5, \quad A=B=1;\\
r_1=r_2=\gamma_1=1, \quad \gamma_2=0.5, \quad \delta_1=1\delta_2=\mu_1=1,
\quad \mu_2=3, \\ m_{21}=30, \quad n_{21}=A=B=1; \\
r_1=r_2=1, \quad \gamma_1=0.5, \quad \gamma_2=\delta_1=1, \quad \delta_2=2, \quad
\mu_1=1, \quad \mu_2=2, \\ m_{21}=1, \quad n_{21}=0.5, \quad A=1, \quad B=3.
\end{eqnarray*}

\section{No Infected Migrations}
In this case the model with $n_{12}=n_{21}=0$ is shown in Figure \ref{fig:graf_interr_migr_2_1}
right and reads
\begin{eqnarray}\label{model_2}
\dot{S_1}=r_1S_1-\gamma_{1}S_1I_1+\delta_1I_1-m_{21}\frac {S_1}{A+I_2+S_2}+m_{12}\frac{S_{2}}{A+S_{1}+I_1}, \\ \nonumber
\dot{I_1}=\gamma_{1}S_1I_1-(\delta_1+\mu_1)I_1 , \\ \nonumber
\dot{S_2}=r_2S_2-\gamma_{2}S_2I_2+\delta_2I_2+m_{21}\frac {S_1}{A+I_2+S_2}-m_{12}\frac{S_{2}}{A+S_{1}+I_1}, \\ \nonumber
\dot{I_2}=\gamma_{2}S_2I_2-(\delta_2+\mu_2)I_2.
\end{eqnarray}

In addition to the origin and ecosystem survival with the endemic disease in both
patches, we find two more points
$$
U=\left(
\frac{\delta_1+\mu_1}{\gamma_1},
r_1\frac{\delta_1+\mu_1}{\gamma_1\mu_1}+\frac{r_2\tilde{S_2}}{\mu_1},
\tilde{S_2}^U,0\right),
\quad
W=\left(
\tilde{S_1}^W,0,
\frac{\delta_2+\mu_2}{\gamma_2}\\
r_2\frac{\delta_2+\mu_2}{\gamma_2\mu_2}+\frac{r_1\tilde{S_1}}{\mu_2}
\right),
$$
where $\tilde{S_2}^U$ and $\tilde{S_1}^W$ solve the equations
\begin{eqnarray}
(r_2\tilde{S_2}A+r_2{\tilde{S_2}}^2+m_{21}\tilde{S_1})\tilde{I_1}=
m_{12}\tilde{S_2}(A+\tilde{S_2})\\  \nonumber
-r_2\tilde{S_2}(A^2+A\tilde{S_1}+A\tilde{S_2}+\tilde{S_1}\tilde{S_2})
-m_{21}\tilde{S_1}(A+\tilde{S_1}),\\
(r_1\tilde{S_1}A+r_1{\tilde{S_1}}^2+m_{12}\tilde{S_2})\tilde{I_1}=
m_{21}\tilde{S_1}(A+\tilde{S_1})\\ \nonumber
-r_1\tilde{S_1}(A^2+A\tilde{S_1}+A\tilde{S_2}+\tilde{S_1}\tilde{S_2})
-m_{12}\tilde{S_2}(A+\tilde{S_2}).
\end{eqnarray}
The first is an intersection problem of the curves
$g(\tilde{S_2})=a_2{\tilde{S_2}}^2+a_1\tilde{S_2}+a_0$,
$f(\tilde{S_2})=b_3{\tilde{S_2}}^3+b_2{\tilde{S_2}}^2+b_1\tilde{S_2}+b_0$, with
\begin{eqnarray}  \label{coef}
a_2=m_{12}-r_2A-r_2\frac{\delta_1+\mu_1}{\gamma_1}, \quad
a_1=Aa_2, \\ \nonumber
a_0=-m_{21}A\frac{\delta_1+\mu_1}{\gamma_1}-m_{21}\frac{{(\delta_1+\mu_1)}^2}{{\gamma_1}^2}<0  \\
b_3=\frac{{r_2}^2}{\mu_1}, \quad
b_2=r_1r_2\frac{\delta_1+\mu_1}{\gamma_1\mu_1}+\frac{{r_2}^2A}{\mu_1}, \\ \nonumber
b_1=r_1r_2A\frac{\delta_1+\mu_1}{\gamma_1\mu_1}+r_2m_{21}\frac{\delta_1+\mu_1}{\gamma_1\mu_1}, \quad
b_0=r_1m_{21}\frac{{(\delta_1+\mu_1)}^2}{{\gamma_1}^2\mu_1}.
\end{eqnarray}
The parabola $g$ has roots ${\tilde{S_2}}^{\pm}=-[a_1{\gamma_1}^2A \pm \sqrt{{a_1}^2{\gamma_1}^4A^2+a_2K}](2a_2 \gamma_1)^{-2}$ and an
intersection for $\tilde S_2\ge 0$ is possible only for $a_2>0$. But the existence
of the intersection is not ensured, since the cubic has positive coefficient,
$b_i\ge 0$, $i=0,...,3$. Similar remarks hold for the
point $W$. Necessary conditions for feasibility are respectively
$$
\tilde{S_2}>{\tilde{S_2}}^{-}=-\frac{a_1{\gamma_1}^2-\sqrt{{a_1}^2{\gamma_1}^4+a_2K}}{2a_2{\gamma_1}^2}, \quad
\tilde{S_1}>{\tilde{S_1}}^{-}=-\frac{c_1\gamma_2A-\sqrt{{c_1}^2{\gamma_2}^2A^2+c_2H}}{2c_2\gamma_2}.
$$

For feasibility of the coexistence equilibrium the necessary conditions (\ref{coex_1}),
(\ref{coex_1bis}) and (\ref{coex_2}) still hold.

The Jacobian of (\ref{model_2}) is
$$
J=
\left[
\begin{array}{cccc}
-\gamma_1I_1-\alpha_1-\beta_2S_2+r_1&-\gamma_1S_1+\delta_1-\beta_2S_2&\alpha_2S_1+\beta_1&\alpha_2S_1 \\
\gamma_1I_1&\gamma_1S_1-\delta_1-\mu_1&0&0\\
\alpha_1+\beta_2S_2&\beta_2S_2& -\gamma_2I_2-\alpha_2S_1+\beta_1+r_2& -\gamma_2S_2+\delta_2-\alpha_2S_1\\
0&0&\gamma_2I_2&\gamma_2S_2-\delta_2-\mu_2
\end{array}
\right]
$$
where
$$
\alpha_1=\frac{m_{21}}{A+I_2+S_2}, \quad \alpha_2=\frac{m_{21}}{{(A+I_2+S_2)}^2},
\quad
\beta_1=\frac{m_{12}}{A+I_1+S_1}, \quad \beta_2=\frac{m_{12}}{{(A+I_1+S_1)}^2}.
$$
The origin is unstable as the eigenvalues are $-\delta_1-\mu_1<0$,
$-\delta_2-\mu_2<0$,
$[k \pm \sqrt{k^2+4A(r_1m_{12}+r_2m_{21}r_2r_2A)}]A^{-1}$, for which one is positive
independently of the sign of $k=r_1A+r_2A-m_{12}-m_{21}$.

At $U$, one eigenvalue is explicit, the other ones are the roots of the cubic
$\sum _{k=0}^3 \lambda^kp_k=0$, for which the Routh-Hurwitz conditions give stability for
\begin{equation}
\label{stab_noinf_caso12}
p_0>0, \quad
p_2>0, \quad
p_2p_1>p_0, \quad
\gamma_2\tilde{S_2}<\delta_2+\mu_2
\end{equation}
where
\begin{eqnarray*}
p_2=\gamma_1\tilde{I_1}+\hat{\alpha_1}+\beta_2-r_1-\gamma_1\tilde{S_1}+\delta_1+\mu_1+\hat{\alpha_2}-\beta_1-r_2, \\ \nonumber
p_1=-(\hat{\alpha_2}+\beta_1)(\hat{\alpha_1}+\beta_2)-\gamma_1\tilde{I_1}(\delta_1-\gamma_1\tilde{S_1}-\beta_2)\\
+(r_1-\gamma_1\tilde{I_1}-\hat{\alpha_1}-\beta_2)
(\gamma_1\tilde{S_1}
-\delta_1-\mu_1)\\ \nonumber
+(r_1-\gamma_1\tilde{I_1}-\hat{\alpha_1}-\beta_2)(r_2-\hat{\alpha_2}+\beta_1)
+(\gamma_1\tilde{S_1}-\delta_1-\mu_1)(r_2-\hat{\alpha_2}+\beta_1), \\
p_0=(r_2-\hat{\alpha_2}+\beta_1)(
+\gamma_1\tilde{I_1}(\delta_1-\gamma_1\tilde{S_1}-\beta_2)\\ \nonumber
-(\hat{\alpha_2}+\beta_1)(\gamma_1\beta_2\tilde{I_1}-(\gamma_1\tilde{S_1}-\delta_1-\mu_1)(\hat{\alpha_1}+\beta_2))\\
-(r_1-\gamma_1\tilde{I_1}-\hat{\alpha_1}-\beta_2)(\gamma_1\tilde{S_1}-\delta_1-\mu_1)
) 
\end{eqnarray*}
A similar situation occurs for $W$, giving stability for
$\gamma_1\tilde{S_1}<\delta_1+\mu_1$ in place of the last one (\ref{stab_noinf_caso12}),
and the same other conditions with $q_i$ in place of $p_i$, where
\begin{eqnarray*}
q_2=\alpha_1+\hat{\beta_2}-r_1+\gamma_2\tilde{I_2}-\alpha_2-\hat{\beta_1}-r_2-\gamma_2\tilde{S_2}+\delta_2+\mu_2 , \\ \nonumber
q_1=(\gamma_2\tilde{S_2}-\delta_2-\mu_2)(r_2-\alpha_1-\hat{\beta_2}+r_1-\gamma_2\tilde{I_2}+\alpha_2+\hat{\beta_1}) \\
    +(r_1-\alpha_1-\hat{\beta_2})(r_2-\gamma_2\tilde{I_2}+\alpha_2+\hat{\beta_1})\\ \nonumber
-(\alpha_1+\hat{\beta_2})(\alpha_2
    +\hat{\beta_1})-\gamma_2\tilde{I_2}( \delta_2-\gamma_2\tilde{S_2}-\alpha_2) , \quad
q_0=(\alpha_2+\hat{\beta_1})(\alpha_1+\hat{\beta_2}))\\
+(\gamma_2\tilde{S_2}-\delta_2-\mu_2)((\alpha_1+\hat{\beta_2}-r_1)(-\gamma_2\tilde{I_2}\\ \nonumber
+\alpha_2+\hat{\beta_1}+r_2)
+\gamma_2\tilde{I_2}((-\alpha_1-\hat{\beta_2}+r_1)(-\gamma_2\tilde{S_2}+\delta_2
    -\alpha_2)-\alpha_2(\alpha_1+\hat{\beta_2})).
\end{eqnarray*}

Numerical simulations reveal that $U$, $W$ and coexistence can be stably achieved, respectively
for the values
\begin{eqnarray*}
r_1=1, \quad r_2=0.5, \quad \gamma_1=0.5, \quad \gamma_2=\delta_1=1, \quad \delta_2=0.2,
\quad \mu_1=1, \\ \mu_2=m_{21}=10, \quad m_{12}=17, \quad A=10;\\
r_1=0.2, \quad r_2=\gamma_1=1, \quad \gamma_2=2, \quad \delta_1=1.8, \quad \delta_2=0.3,
\quad
\mu_1=1, \\ \mu_2=6, \quad m_{21}=8.8, \quad m_{12}=4, \quad A=10;\\
r_1=r_2=\gamma_1=\gamma_2=1, \quad \delta_1=2, \quad \delta_2=0.5, \quad
\mu_1=2\mu_2=2, \\ m_{12}=3, \quad m_{21}=1, \quad A=10.
\end{eqnarray*}

\section{Interpretation}

The generic equilibria share common properties in both models. Namely,
instability of the origin means ecosystem permanence. The systems allow also
survival of the population in both patches, with endemic disease.

Within the unidirectional migration model, the patch from which migrations occur could become
completely depleted, or else it may be populated, but disease-free. 
As these are mutually exclusive equilibria, it is possible that bistability phenomena
arise as for \cite{ABMV}.
The basins of attraction of the equilibria could be determined
using the algorithms being developed, \cite{CCDM}.

For the model in which infected do not migrate, specific equilibria are the situations
in which either patch becomes epidemic-free.
These considerations could be very useful in practical situations for disease
eradication in some environments.

\end{document}